\documentclass{elsart}
\usepackage{amssymb,amsfonts}
\newcommand{\trace}{\mathop{\rm Tr}\nolimits}

\newcommand{\twomat}[4]{\left(\begin{array}{cc}#1&#2\\#3&#4\end{array}\right)}
\newcommand{\twovec}[2]{\left(\begin{array}{c}#1\\#2\end{array}\right)}

\newcommand{\cP}{{\mathcal P}} % \cal not known

\newcommand{\R}{{\mathbb{R}}}
\newcommand{\N}{{\mathbb{N}}}
\newcommand{\id}{\mathbb{I}}

\newcommand{\be}{\begin{equation}}
\newcommand{\ee}{\end{equation}}
\newcommand{\bea}{\begin{eqnarray}}
\newcommand{\eea}{\end{eqnarray}}
\newcommand{\beas}{\begin{eqnarray*}}
\newcommand{\eeas}{\end{eqnarray*}}
\newtheorem{definition}{Definition}
\newtheorem{theorem}{Theorem}
\newtheorem{lemma}{Lemma}
\newtheorem{corollary}{Corollary}
\newtheorem{conjecture}{Conjecture}

\newtheorem{proposition}{Proposition}
\newcount\minute
\newcount\hour
\def\currenttime{%
    \minute\time
    \hour\minute
    \divide\hour60
    \the\hour:\multiply\hour60\advance\minute-\hour\the\minute}
\begin{document}
\begin{frontmatter}
%\title{Complementary McCarthy inequalities}
\title{Trace inequalities for completely monotone functions and Bernstein functions}
\author{Koenraad M.R. Audenaert}
\address{
Department of Mathematics, Royal Holloway, University of London, Egham TW20 0EX, United Kingdom}
\ead{koenraad.audenaert@rhul.ac.uk}
\date{\today, \currenttime}
%------------------------------------------------------------------ ABSTRACT
\begin{abstract}
We prove a matrix trace inequality for completely monotone functions and for Bernstein functions.
As special cases we obtain non-trivial trace inequalities for the power function $x\mapsto x^q$,
which for certain values of $q$
complement McCarthy's trace inequality and for others strenghten it.
\end{abstract}

\end{frontmatter}
%%%%%%%%%%%%%%%%%%%%%%%%%%%%%%%%%%%%%%%%%%%%%%%%%%%%%%%%%%%%%%%%%%%%%%%%%%%%%%%%%%%%%%%%%%%%%
\section{Introduction\label{sec1}}
Completely monotone functions play an important role in many branches of applied mathematics, and probability
theory. They are defined as the Laplace transforms of measures on the half-line $[0,\infty)$.
Closely related to these functions are the so-called Bernstein functions, which are the primitives of the
positive completely monotone functions. 
Bernstein functions are important in probability theory as well, appearing
for example in the study of random Markov processes.

In Section \ref{sec:main} of this paper we prove a matrix trace inequality
for completely monotone functions, Bernstein functions, and
primitive functions thereof. Our interest in these classes of functions
does not stem from the applications just mentioned but from the fact
that they contain the fractional power function $x\mapsto x^q$,
for various ranges of $q$. 

By specialising our trace inequality to the fractional power function
we obtain in Section \ref{sec:power} a number of non-trivial trace inequalities
related to McCarthy's trace inequality \cite{cp}. The latter states that the
matrix function $A\mapsto\trace A^q$ is subadditive on the set of
positive semidefinite matrices for $0<q\le1$, and superadditive for $q\ge1$:
\bea
\trace(A+B)^q&\le&\trace A^q+\trace B^q,\qquad 0<q\le 1 \\
\trace(A+B)^q&\ge&\trace A^q+\trace B^q,\qquad q\ge 1.
\eea
The inequalities we obtain are complementary to McCarthy's in particular regions for $q$,
and are strengthenings in others. 

As a further application of these inequalities we obtain in Section \ref{sec:norm} 
a simple proof of a norm inequality for partitioned positive
semidefinite matrices that was first proven in \cite{ka} by other means.
%%%%%%%%%%%%%%%%%%%%%%%%%%%%%%%%%%%%%%%%%%%%%%%%%%%%%%%%%%%%%%%%%%%%%%%%%%%%%%%%%%%%%%%%%%%%%%%%%%%%%%%%%%
%\section{Preliminaries\label{sec2}}

%%%%%%%%%%%%%%%%%%%%%%%%%%%%%%%%%%%%%%%%%%%%%%%%%%%%%%%%%%%%%%%%%%%%%%%%%%%%%%%%%%%%%%%%%%%%%%%%%%%%%%%%%%%%%%%%%%%

\section{Completely monotone functions and Bernstein functions}
In this section, we collect a number of definitions and theorems about completely monotone functions
and Bernstein functions that will be needed later on.
For an in-depth treatment, including proofs and applications, we refer to the excellent monograph \cite{bernstein}.

\subsection{Completely monotone functions}
\begin{definition}
A function $f:(0,\infty)\mapsto\R$ is \textit{completely monotone} if it is infinitely differentiable,
non-negative, and $(-1)^n f^{(n)}(x)\ge 0$ for $n=1,2,\ldots$ and $x>0$.
\end{definition}

An integral representation of completely monotone functions is provided by Bernstein's theorem:
\begin{theorem}[Bernstein]
A function $f:(0,\infty)\mapsto\R$ is completely monotone if and only if it is the Laplace transform
of a positive measure $\mu$ on $[0,\infty)$, i.e.
\be
f(x) = \int_{[0,\infty)}\exp(-x t) \mu(dt)
= a+\int_{(0,\infty)}\exp(-x t) \mu(dt), \label{eq:cm}
\ee
where $a$ is given by $a=f(0^+)$.
\end{theorem}

In addition, we will also define the \textit{bare} completely monotone functions as
those completely monotone functions for which $a=f(0^+)=0$.
We will denote the class of bare completely monotone functions by $\mathsf{CM}0$.

An important class of completely monotone functions are the negative power functions
$x\mapsto x^{q}$, $q<0$.
That these functions are completely monotone follows from the integral representation
\be
x^{q} = \frac{1}{\Gamma(-q)}\;\int_0^\infty \exp(-xt)t^{-q-1}\,dt,\qquad q<0.\label{eq:xqcm}
\ee
This fits the representation of Bernstein's theorem via $\mu(dt)=t^{-q-1}dt/\Gamma(-q)$.

\begin{lemma}\label{lem:CM0}
Any function $f\in\mathsf{CM}0$ is convex, monotonically decreasing and non-negative.
\end{lemma}
\textit{Proof.}
Obvious from the integral representation
$$
f(x)=\int_{(0,\infty)}\exp(-x t) \mu(dt),
$$
since $\exp(-x)$ is convex, monotonically decreasing and non-negative.
\qed
%%%%%%%%%%%%%%%%%%%%%%%%%%%%%%%%%%%%
\subsection{Bernstein functions}
\begin{definition}
A function $f:(0,\infty)\mapsto\R$ is a \textit{Bernstein function} if it is infinitely differentiable,
non-negative, and $(-1)^n f^{(n)}(x)\le 0$ for $n=1,2,\ldots$ and $x>0$.
\end{definition}

Again, this class of functions admits an integral representation.
\begin{theorem}\label{th:LK}
A function $f:(0,\infty)\mapsto\R$ is a Bernstein function if and only if
there exist $a,b\ge0$ and a positive measure $\mu(dt)$ on $(0,\infty)$
such that
\be
f(x) = ax+b + \int_{(0,\infty)}\left(1-\exp(-x t)\right) \mu(dt). \label{eq:bf}
\ee
\end{theorem}
In the probability theory literature, this representation is known as the L\'evy-Khintchine representation.
The constants $a$ and $b$ are given by the limits
$a=\lim_{x\to\infty}f(x)/x$ and $b=f(0^+)$.
In addition, we define the \textit{bare} Bernstein functions
as those Bernstein functions for which $a=b=0$, and denote this class by $\mathsf{BF}0$.

The kernel function $1-\exp(-x t)$ has leading order degree 1 for $x$ tending to 0, and
degree 0 for $x$ tending to $\infty$.
Therefore, in order for the integral in (\ref{eq:bf}) to converge the measure $\mu$
must satisfy in integrability condition, such as the following one:
$$
\int_{(0,\infty)}\min(1,t)\mu(dt)<\infty
% \mbox{ or }
%\int_{(0,\infty)}\frac{t}{1+t}\mu(dt)<\infty.
$$

It is easy to see that the derivative of every Bernstein function is completely monotone.
Indeed, representation (\ref{eq:bf}) can be obtained from (\ref{eq:cm}) by simple integration.
However, not every completely monotone function is the derivative of a Bernstein function,
because of the extra positivity requirement for the latter.
The function $x\mapsto x^q$, $q\le -1$, for example, is completely monotone 
but its primitive $x^{q+1}/(q+1)$
is negative for all $x>0$ and therefore not a Bernstein function.

An important subclass of the Bernstein functions are the non-negative operator monotone functions.
They contain the fractional power functions
$x\mapsto x^q$, for $0<q<1$, as can be seen from the integral
representation
\be
x^q = \frac{q}{\Gamma(1-q)}\int_0^\infty(1-\exp(-tx))t^{-q-1}\,dt,\qquad 0<q<1. \label{eq:xqbf}
\ee

\begin{lemma}\label{lem:BF0}
Any function $f\in\mathsf{BF}0$ is concave, monotonically increasing, non-negative and $f(0)=0$.
\end{lemma}
\textit{Proof.}
Obvious from the integral representation
$$
f(x)=\int_{(0,\infty)}(1-\exp(-x t)) \mu(dt),
$$
since $1-\exp(-x)$ is concave, monotonically increasing and non-negative, and $1-\exp(0)=0$.
\qed

\subsection{Integrals of Bernstein functions}
We will also consider functions whose first (second) derivative is a bare Bernstein function.
More precisely, given any bare Bernstein function $f\in\mathsf{BF}0$, we consider the functions
$$
g(y) = \int_0^y dx\, f(x) \mbox{ and } h(z)=\int_0^z dy\, g(y).
$$
From representation (\ref{eq:bf}) it follows that the function $g$
is represented by
\be
g(x) = \int_{(0,\infty)}\left(\exp(-x t)-(1-xt)\right)\,\frac{1}{t} \mu(dt), \label{eq:bfp}
\ee
with $\mu(dt)$ the measure appearing in Theorem \ref{th:LK}.
The class of these functions will be denoted by $\mathsf{BF}1$.
For $x$ tending to 0, the leading order of the kernel function $\exp(-x t)-(1-xt)$ is of degree $2$,
while for $x$ tending to $\infty$ it is of degree $1$.
Convergence of the integral is therefore not affected by the extra factor $1/t$.

Likewise, the function $h$ is represented by
\be
h(x) =
\int_{(0,\infty)}\left((1-xt+x^2t^2/2)-\exp(-x t)\right)\,\frac{1}{t^2} \mu(dt). \label{eq:bfpp}
\ee
The class of these functions will be denoted by $\mathsf{BF}2$.
Here, the kernel function $\exp(-x t)-(1-xt+x^2t^2/2)$ has leading order degree 3 for $x$ tending to 0,
so that convergence is again not affected by the factor $1/t^2$.

Continuing in this way, we can inductively define the classes $\mathsf{BF}k$, $k\in\N$,
as the classes of $k$-fold integrals of bare Bernstein functions.
That is, $f\in\mathsf{BF}k$ if and only if there is a function $g\in\mathsf{BF}(k-1)$
such that
\be
f(x) = \int_0^x g(t)dt. \label{eq:bfkdef}
\ee
It is easy to see that such functions have the integral representation
\be
f(x) = \int_{(0,\infty)}
(-1)^{k+1} \left(\exp(-x t)-\sum_{j=0}^k(-xt)^j/j!\right)\,\frac{1}{t^k} \mu(dt). \label{eq:bfk}
\ee
For $x$ tending to 0, the leading order of the kernel function is of degree $k+1$,
while for $x$ tending to $\infty$ it is of degree $k$.

The function $x\mapsto x^q$ is
%in $\mathsf{BF1}$ for $1< q< 2$, in $\mathsf{BF2}$ for $2< q< 3$, and in general
in $\mathsf{BF}k$ for $k<q<k+1$.

\begin{lemma}\label{lem:BFk}
Any function $f\in\mathsf{BF}k$, $k\ge1$, is convex, monotonically increasing, non-negative, and $f(0)=0$.
\end{lemma}
\textit{Proof.}
For $k=1$, this is obvious from integral representation (\ref{eq:bfp}),
since the function $x\mapsto\exp(-x)-1+x$ is convex, monotonically increasing, non-negative and $\exp(0)-1+0=0$.

For $k>1$, this follows inductively from the defining integral (\ref{eq:bfkdef}), from which we get
$f\ge0$, $f(0)=0$, $f'=g$ and $f''=g'$, for $g\in\mathsf{BF}(k-1)$.
By the induction hypothesis, $g$ is non-negative and increasing, hence
$f'\ge0$ and $f''\ge0$.
\qed

%%%%%%%%%%%%%%%%%%%%%%%%%%%%%%%%%%%%%%%%%%%%
\section{Main Results\label{sec:main}}
As stated in the Introduction, we will exploit the integral representations of functions in $\mathsf{CM}0$
and $\mathsf{BF}k$ to extend inequalities for the exponential function to those classes of functions.
\subsection{Scalar inequalities\label{sec:main:scalar}}
In this section, we restrict ourselves to the scalar case, leaving the matrix case for the next section.
%%%%%%%%%%%%%%%%%%%%%%%%%%%%%%%%%%%%%%%%%%%%%%%%%%%%%%%%%%%%%%%%%%%%%%%%%%%%%

The following lemma concerns (scalar) sub- and superadditivity.
Recall that a real-valued function $g$ is subadditive on $I$ if and only if
$\forall x,y\in I: g(x+y)\le g(x)+g(y)$; it is superadditive on $I$ if and only if
$\forall x,y\in I: g(x+y)\ge g(x)+g(y)$.
\begin{lemma}\label{lem:sub1}
Let $g$ be a function $g:[0,\infty)\to\R$.
If $g\in\mathsf{CM}0\cup\mathsf{BF}0$ then $g$ is subadditive on $[0,\infty)$.
If $g\in\mathsf{BF}k$, $k\ge1$, then $g$ is superadditive on $[0,\infty)$.
\end{lemma}
\textit{Proof.}
For $x,y\ge0$, we have $e^x+e^y\ge 2$, hence
$e^{-x-y} \le 2e^{-x-y}\le e^{-y}+e^{-x}$. Thus, the function $e^{-xt}$ is subadditive for all $t\ge0$.
Therefore, all functions in $\mathsf{CM}0$ are subadditive too.

The same is true for functions in $\mathsf{BF}0$, as can be seen from subadditivity of the function
$1-e^{-x}$.
The latter follows from positivity of $(1-e^{-x})(1-e^{-y})$ for $x,y\ge0$.

Superadditivity of functions in $\mathsf{BF}1$ follows from superadditivity of the function
$x\mapsto e^{-x}-1+x$, which in turn follows from subadditivity of $1-e^{-x}$ and additivity of $x\mapsto x$.

Superadditivity of functions in $\mathsf{BF}k$ for $k>1$ follows inductively from superadditivity of
functions in $\mathsf{BF}(k-1)$.
By definition, any function $h\in\mathsf{BF}k$ is given by the integral
$h(x)=\int_0^x dt g(t)$ of a function $g\in\mathsf{BF}(k-1)$.
Therefore,
\beas
h(x+y)-h(x)-h(y)&=&\int_x^{x+y} dt g(t) - \int_0^{y} dt g(t) \\
&=& \int_0^{y} dt (g(t+x)-g(t)) \\
&\ge& \int_0^{y} dt g(x) = y g(x)\ge0.
\eeas
In the last line we exploited superadditivity of $g$ in the form $g(t+x)-g(t)\ge g(x)$.
\qed

%%%%%%%%%%%%%%%%%%%%%%%%%%%%%%%%%%%%%%%%%%%%%%%%%%%%%%%%%%%%%%%%%%%%%%%%%%%%%

To obtain inequalities that complement the subadditivity (superadditivity) inequalities of the previous
lemma, we need a not very well-known property of the exponential function.
This property relies on the so-called geometrical concavity of the function $1-\exp(-x)$,
a concept that can be traced back to Montel \cite{montel}.
\begin{definition}
A function $f(x): \R^+\to\R^+$ is \textit{geometrically concave} iff for all $x,y\ge0$,
$$
f(\sqrt{xy})\ge\sqrt{f(x)f(y)}.
$$
\end{definition}
\begin{lemma}\label{lem:geoconc}
The function $f(x)=1-\exp(-x)$ is geometrically concave.
%, and satisfies $0\le f'(x)\le f(x)/x$ for all $x>0$.
\end{lemma}
\textit{Proof.}
Geometrical concavity of $f$ is equivalent to concavity of
$g(x)=\log(1-\exp(-\exp x))$.
The second order derivative of $g$ is
$$
g''(x) = \frac{\exp x}{(\exp\exp x-1)^2}(\exp\exp x-\exp(x+\exp x)-1).
$$
The factor that determines the sign is clearly $\exp\exp x-\exp(x+\exp x)-1$, which is non-positive.
Indeed, substituting $a=\exp x$, and noting that $\exp(-a)\ge 1-a$, yields
\beas
\exp\exp x-\exp(x+\exp x)-1 &=& \exp a-a\exp a-1 \\
&=& (1-a)\exp a-1 \\
&\le& \exp(-a)\exp a-1 = 0.
\eeas
%
%Obviously, $f$ is non-decreasing, hence $f'\ge0$.
%
%Finally, $f'(x) = \exp(-x)$ and $f(x)/x = (1-\exp(-x))/x$. Since $\log(1+x)\le x$ for all $x\ge0$,
% $\exp(-x)\le 1/(1+x)$. Thus, $(1+x)\exp(-x)\le 1$, or $\exp(-x) \le (1-\exp(-x))/x$. This shows
%that $f'(x)\le f(x)/x$ for all $x>0$.
\qed

This property of the function $1-\exp(x)$  translates to a property of the
exponential function. The connection to subadditivity (superadditivity) is immediate.
\begin{lemma}\label{lem:cmf1}
The inequality
\be
g(a+b) - g(a) - g(b) \le g(2\sqrt{ab}) - 2g(\sqrt{ab})
\label{eq:cmf1a}
\ee
holds for any $a,b\ge0$ when $g(x)=\exp(-x)$.
If $g(x)$ is a quadratic polynomial then it holds with equality.
\end{lemma}
\textit{Proof.}
That equality holds for quadratic polynomials is immediate.

Geometrical concavity of $f(x)=1-\exp(-x)$ amounts to the inequality
$$
(1-\exp(-a))(1-\exp(-b)) \le (1-\exp(-\sqrt{ab}))^2.
$$
Hence, for all $a,b\ge0$,
\beas
\exp(-(a+b)) - \exp(-a) - \exp(-b) &\le& \exp(-2\sqrt{ab}) - 2\exp(-\sqrt{ab}).
\eeas
\qed

This inequality can be extended to completely monotone functions and Bernstein functions,
using their integral representations.
\begin{theorem}\label{th:cmf1}
Let $g$ be a function $g: [0,\infty)\to\R$.
For $a,b\ge0$,
\be
g(a+b) - g(a) - g(b) \le g(2\sqrt{ab}) - 2g(\sqrt{ab}),\label{eq:cmf1}
\ee
holds if $g\in\mathsf{CM}0\cup\mathsf{BF}1$.
The inequality holds in the reversed sense if $g\in\mathsf{BF}0\cup\mathsf{BF}2$.
Equality holds when $g$ is a quadratic polynomial.
\end{theorem}
\textit{Proof.}
This follows immediately from Lemma \ref{lem:cmf1} and from the integral representations of functions in
$\mathsf{CM}0$, $\mathsf{BF}0$, $\mathsf{BF}1$ and $\mathsf{BF}2$.
The sign with which $\exp(-tx)$ occurs in these representations determines whether the inequality holds
in the stated sense or in the reversed sense.
\qed

It will be shown below that this inequality does not hold for functions in $\mathsf{BF}k$, $k>2$.

The results from Lemma \ref{lem:sub1} and Theorem \ref{th:cmf1} can be summarised by
the following inequalities:
$$
\begin{array}{cccccccl}
&&g(x+y)-g(x)-g(y)&\le&g(2\sqrt{xy})-2g(\sqrt{xy})&\le&0, & \mbox{ } g\in\mathsf{CM}0 \\
0&\ge&g(x+y)-g(x)-g(y)&\ge&g(2\sqrt{xy})-2g(\sqrt{xy})&&, & \mbox{ } g\in\mathsf{BF}0 \\
0&\le&g(x+y)-g(x)-g(y)&\le&g(2\sqrt{xy})-2g(\sqrt{xy})&&, & \mbox{ } g\in\mathsf{BF}1 \\
&&g(x+y)-g(x)-g(y)&\ge&g(2\sqrt{xy})-2g(\sqrt{xy})&\ge&0, & \mbox{ } g\in\mathsf{BF}2.
\end{array}
$$
Thus, Theorem \ref{th:cmf1} is a stronger statement than subadditivity (superadditivity)
for $g\in \mathsf{CM}0$ ($g\in\mathsf{BF}2$), while for $g\in \mathsf{BF}0$ ($g\in\mathsf{BF}1$) it
provides a complementary inequality to subadditivity (superadditivity).

%%%%%%%%%%%%%%%%%%%%%%%%%%%%%%%%%%%%%%%%%%%%%%%%%%%%%%%%%%%%%%%%%%%%%%%%%%%%%
\subsection{A matrix trace inequality\label{sec:main:matrix}}
Theorem \ref{th:cmf1} is easy to extend to the matrix case via a simple application of
the Golden-Thompson theorem, yielding our main trace inequality:
\begin{theorem}\label{th:cmf2}
Let $A$ and $B$ be $d$-dimensional positive semidefinite matrices, with spectral decompositions
$A=\sum_k a_k A_k$ and $B=\sum_k b_k B_k$, where $a_k, b_k\ge0$ and
$\{A_k\}$ and $\{B_k\}$ are two complete sets of mutually orthogonal projectors.
Let $g$ be a function, $g: [0,\infty)\to\R$.
The inequality
$$
\trace(g(A+B) - g(A) - g(B)) \le \sum_{k,l}(g(2\sqrt{a_kb_l}) - 2g(\sqrt{a_kb_l}))\trace A_kB_l
$$
holds if $g\in\mathsf{BF}1$ and (for $A,B>0$) if $g\in\mathsf{CM}0$.

The inequality holds in the reversed sense if $g\in\mathsf{BF}0\cup\mathsf{BF}2$.

Equality holds if $g$ is a quadratic polynomial.
\end{theorem}
\textit{Proof.}
It is easy to check that the inequality reduces to an equality for $g(x)=1$, $g(x)=x$ and $g(x)=x^2$.
For $g(x)=1$, the LHS is $-\trace\id$, and the RHS is
$-\sum_{k,l}\trace A_kB_l = -\trace(\sum_k A_k)(\sum_l B_l) = -\trace\id$, due to completeness of the sets
$\{A_k\}$ and $\{B_k\}$.

For $g(x)=x$, the LHS and RHS are both 0, and for $g(x)=x^2$ the LHS is $2\trace AB$
and the RHS is $\sum_{k,l}2a_kb_l\trace A_kB_l = 2\trace(\sum_k a_k A_k)(\sum_l b_l B_l) = 2\trace AB$.

To prove the main statement of the theorem, we look again at the exponential function.
The Golden-Thompson theorem states
$\trace\exp(A+B)\le \trace\exp A\exp B$, for any two Hermitian matrices $A$ and $B$.
In particular, we have, for any $t$,
\beas
\trace\exp(-(A+B)t)&\le& \trace\exp(-At)\exp(-Bt)\\
 &=& \sum_{k,l} \exp(-a_kt)\exp(-b_lt)\trace A_kB_l\\
 &=& \sum_{k,l} \exp(-(a_k+b_l)t)\trace A_kB_l.
\eeas
Also, for any function $g$,
\beas
\trace g(A) &=& \sum_{k} g(a_k)\trace A_k = \sum_{k,l} g(a_k)\trace A_kB_l,\\
\trace g(B) &=& \sum_{l} g(b_l)\trace B_l = \sum_{k,l} g(b_l)\trace A_kB_l.
\eeas
Therefore, for $g(x)=\exp(-xt)$,
$$
\trace (g(A+B)-g(A)-g(B)) \le \sum_{k,l}(g(a_k+b_l)-g(a_k)-g(b_l))\trace A_kB_l.
$$

Using the same reasoning as in the proof of Theorem \ref{th:cmf1}, we then find that this inequality
holds for all $g\in\mathsf{BF}1$ and (for $A,B>0$) $g\in\mathsf{CM}0$,
and in the reversed sense for all $g\in\mathsf{BF}0\cup\mathsf{BF}2$.

Combining this with the scalar inequality of Theorem \ref{th:cmf1} applied to $g(a_k+b_l)-g(a_k)-g(b_l)$
yields the stated inequalities.
\qed

For completeness, we also state the extension of Lemma \ref{lem:sub1}
to the matrix case.
Let $\cP$ denote the set of positive semidefinite matrices.
\begin{lemma}\label{lem:sub2}
Let $g$ be a function $g:[0,\infty)\to\R$, and extended to $\cP$ in the usual way.
If $g\in\mathsf{CM}0\cup\mathsf{BF}0$ then the function $A\mapsto\trace g(A)$ is subadditive on $\cP$,
i.e.\ for all $A,B\ge0$,
$$
\trace g(A+B) \le \trace(g(A)+g(B)).
$$
If $g\in\mathsf{BF}k$, $k\ge1$, then $A\mapsto\trace g(A)$ is superadditive on $\cP$, i.e.\ for all $A,B\ge0$,
$$
\trace g(A+B) \ge \trace(g(A)+g(B)).
$$
\end{lemma}
\textit{Proof.}
To show the statement for $g\in\mathsf{CM}0$, we only need to show it for $g(x)=\exp(-x)$, i.e.\
that $\trace e^{-A-B}\le \trace e^{-A}+e^{-B}$.
By the Golden-Thompson inequality, we have $\trace e^{-A-B}\le \trace e^{-A}e^{-B}$.
Since $B\ge0$, we also have $e^{-B}\le\id$.
Thus, $\trace e^{-A}e^{-B}\le \trace e^{-A} \le \trace e^{-A} +e^{-B}$,
so that indeed $\trace e^{-A-B}\le \trace e^{-A}+e^{-B}$.

Next, to cover the case $g\in\mathsf{BF}0$, we just note that $\trace g(A+B) \le \trace(g(A)+g(B))$
is a special case of Bourin and Uchiyama's norm subadditivity inequality \cite{BU}. Indeed, by Lemma \ref{lem:BF0},
functions in $\mathsf{BF}0$ satisfy the conditions of their Theorem.

Likewise, to cover the case $g\in\mathsf{BF}k$, $k\ge1$, we note that $\trace g(A+B) \ge \trace(g(A)+g(B))$
is a special case of Kosem's norm superadditivity inequality \cite{kosem}. By Lemma \ref{lem:BFk},
functions in $\mathsf{BF}k$, $k\ge1$, satisfy the conditions of his Theorem.
\qed

%%%%%%%%%%%%%%%%%%%%%%%%%%%%%%%%%%%%%%%%%%%%%%%%%%%%%%%%%%%%%%%%%%%%
\section{Inequalities for the power function\label{sec:power}}
The inequality of Theorem \ref{th:cmf2} achieves its most elegant form when $g(x)$ is the
fractional power function $x\mapsto x^q$.
\begin{corollary}\label{cor:abq}
For $A,B\ge0$ and $0<q\le1$ or $2\le q\le 3$,
\be
\trace(A+B)^q-\trace(A^q+B^q) \ge (2^q-2)\trace A^{q/2}B^{q/2}.\label{eq:abq}
\ee
For $q<0$ (in which case we require $A,B>0$) and $1\le q\le 2$, the inequality holds in the reversed sense.
\end{corollary}
\textit{Proof.}
For $g(x)=x^q$, the RHS of the inequality of Theorem \ref{th:cmf2} simplifies to
\beas
\lefteqn{\sum_{k,l}\left(\left(2\sqrt{a_kb_l}\right)^q - 2\left(\sqrt{a_kb_l}\right)^q\right)\trace A_kB_l} \\
&=& (2^q-2) \sum_{k,l}a_k^{q/2} b_l^{q/2}\trace A_kB_l \\
&=& (2^q-2)\trace A^{q/2} B^{q/2}.
\eeas
The corollary then follows by recalling that
$g(x)=x^q$ is in $\mathsf{CM}0$ for $q<0$, and in $\mathsf{BF}k$ for $k<q<k+1$.
\qed

For $q>3$ the inequality does no longer hold in general.
Indeed, for scalars (e.g. $A=1$, $B=2$) the inequality holds in the stated sense
while, for example, with the choice $A=\twomat{1}{0}{0}{0}$, $B=\frac{1}{2}\twomat{1}{1}{1}{1}$, the inequality
holds in the reversed sense for all $q>3$; the LHS is $(1+\sqrt{2}/2)^q+(1-\sqrt{2}/2)^q-2$ and the RHS
$(2^q-2)/2$.
More generally, this shows that the inequality of Theorem \ref{th:cmf2} does not hold for functions
in $\mathsf{BF}k$ for $k>2$.

As in the scalar case,
one sees that the inequality of Corollary \ref{cor:abq}
is stronger than McCarthy's for $2\le q\le 3$. The corollary also implies that subadditivity holds for
$q<0$ too.
On the other hand, for other parameter ranges the corollary complements McCarthy's inequalities
by providing a lower bound on $\trace(A+B)^q-\trace A^q-\trace B^q$ for $0<q\le1$
and an upper bound for $1\le q\le2$.

By replacing $A$ and $B$ by $A^{1/q}$ and $B^{1/q}$, and $q$ by $1/p$,
Corollary \ref{cor:abq} can be reformulated as an inequality for the $p$-power means \cite{bhagwat}:
\begin{corollary}
For $A,B\ge0$ and $p\ge1$,
$$
\trace\left(\frac{A^p+B^p}{2}\right)^{1/p} \ge  2^{1-1/p}\trace\frac{A+B}{2} + (1-2^{1-1/p})\trace A^{1/2}B^{1/2}.
$$
\end{corollary}

\bigskip

Applying the Araki-Lieb-Thirring inequality, we obtain a closely related trace inequality, where the expression
$\trace A^{q/2}B^{q/2}$ is replaced by $\trace(A^{1/2}BA^{1/2})^{q/2}$.
\begin{corollary}\label{cor:FALTq}\label{cor:abq2}
For $A,B\ge0$ and $0<q\le1$ or $2\le q\le 3$,
\be
\trace(A+B)^q-\trace(A^q+B^q) \ge (2^q-2)\trace(A^{1/2}BA^{1/2})^{q/2}.\label{eq:FALTq}
\ee
For $q\le -2$ (in which case we require $A,B>0$) and $1\le q\le 2$, the inequality holds in the reversed sense.
\end{corollary}
\textit{Proof.}
As always, we require $A,B>0$ for negative $q$.

By the Araki-Lieb-Thirring inequality for $A,B\ge0$, we have,
for $0<q\le2$ and for $-2\le q<0$,
$$
\trace A^{q/2}B^{q/2} \le \trace (A^{1/2}BA^{1/2})^{q/2},
$$
while for $q\ge2$ and for $q\le -2$ the reversed inequality holds;
In addition, $2^q-2$ is positive for $q>1$ and negative for $q<1$.
Thus, the corollary follows from Corollary \ref{cor:abq}
for all values of $q$ for which the latter holds (i.e.\ $q\le 3$),
except for $-2<q<0$.
\qed

In contrast to inequality (\ref{eq:abq}), which does not hold for $q>3$,
we have numerical evidence in support of the following conjecture:
\begin{conjecture}
Inequality (\ref{eq:FALTq}) also holds for $q>3$ and, in the reversed sense, for $-2<q<0$.
\end{conjecture}
Additional evidence is given by:
\begin{proposition}
Inequality (\ref{eq:FALTq}) holds for $q=4$.
\end{proposition}
\textit{Proof.}
For $q=4$, the left-hand side of (\ref{eq:FALTq}) is
$$
\trace(A+B)^4-\trace A^4-\trace B^4 = 4\trace (A^3B+A^2B^2+AB^3)+2\trace(AB)^2,
$$
while the right-hand side is
$$
(2^4-2)\trace(A^{1/2}BA^{1/2})^2 = 12\trace(AB)^2.
$$

Now, $\trace(A^3B+AB^3)\ge 2\trace A^2B^2$.
This follows from the scalar inequality $x^3y+xy^3=xy(x^2+y^2)\ge xy(2xy)=2x^2y^2$
applied to the coefficients in the spectral decomposition
of $\trace(A^3B+AB^3)=\sum_{j,k} (a_j^3b_k +a_jb_k^3)\trace A_jB_k$
(in terms of the spectral decompositions $A=\sum_j a_j A_j$ and $B=\sum_k b_k B_k$).
Note that $\trace A_jB_k\ge0$.

Also, by the Lieb-Thirring inequality, $\trace A^2B^2\ge \trace(AB)^2$.
Thus, indeed we find
$$
4\trace (A^3B+A^2B^2+AB^3)+2\trace(AB)^2 \ge 12\trace(AB)^2
$$
\qed
%%%%%%%%%%%%%%%%%%%%%%%%%%%%%%%%%%%%%%%%%%%%%%%%%%%%%%%%%%%%%%%%%%%%%%%%%%%%%%%%%%%%%%%%%%%%%%%%%%%%%%%%%%
\section{A new proof for a norm compression inequality\label{sec:norm}}
In this section, we consider an inequality that relates the Schatten $q$-norm,
$||X||_q := (\trace|X|^q)^{1/q}$, of a partitioned positive semidefinite
matrix to the Schatten $q$-norms of its blocks. In particular, we compare
it to the $q$-norm of the matrix that is obtained by replacing each block by its $q$-norm.
An inequality of this type is sometimes called a \textit{norm compression inequality}.
The specific inequality presented here has first appeared in our \cite{ka}, but had
a long and intricate proof, and only for the case $1\le q\le 2$.
Here we show how to prove it in a very simple way, using Corollary \ref{cor:abq2} from the previous section.
Moreover, the proof given here extends the result to include the parameter range $0<q<1$ and $2<q\le3$.

In the following, we consider the positive semidefinite block matrix
$$
A=\twomat{B}{C^*}{C}{D},
$$
where $B$ and $D$ are square blocks, and define the block norms
$\beta=||B||_q$, $\delta=||D||_q$ and $\gamma=||C||_q$.

An equivalent form of the inequalities in Corollary \ref{cor:abq2} is:
\begin{corollary}\label{cor:abq3}
Let $D>0$ and let $C$ be any matrix.
For $q\le -2$ and for $1\le q\le 2$,
$$
\trace\twomat{C^*D^{-1}C}{C^*}{C}{D}^q-\trace\twomat{C^*D^{-1}C}{0}{0}{D}^q
\le (2^q-2)\trace |C|^q.
$$
For $0<q\le 1$ and for $2\le q\le 3$ the reversed inequality holds.
\end{corollary}
\textit{Proof.} In the inequalities of Corollary \ref{cor:abq2}
set $A=D^{-1/2}CC^*D^{-1/2}$ and $B=D$.
Then
$$
A+B = (D^{-1/2}C \,\,\, D^{1/2})\twovec{C^*D^{-1/2}}{D^{1/2}},
$$
which has the same non-zero eigenvalues as the block matrix
$$
Z=\twovec{C^*D^{-1/2}}{D^{1/2}}(D^{-1/2}C \,\,\, D^{1/2})
=\twomat{C^*D^{-1}C}{C^*}{C}{D}.
$$
Therefore $\trace (A+B)^q$ is equal to $\trace Z^q$.

Furthermore, $A=D^{-1/2}CC^*D^{-1/2}$ and $C^*D^{-1/2}D^{-1/2}C = C^*D^{-1}C$ are unitarily
equivalent, so that $\trace A^q=\trace(C^*D^{-1}C)^q$.

Finally, because $A=D^{-1/2}CC^*D^{-1/2}$ there exists a unitary matrix $U$ such that
$C=D^{1/2}A^{1/2}U=B^{1/2}A^{1/2}U$.
Thus, $(A^{1/2}BA^{1/2})^{1/2}=U(C^*C)^{1/2}U^* = U|C|U^*$,
whence
$$
\trace(A^{1/2}BA^{1/2})^{q/2}=\trace|C|^q.
$$

Substituting everything in the inequality of Corollary \ref{cor:abq2} yields the stated inequality.

Conversely, the inequality of Corollary \ref{cor:abq2} is obtained from the stated inequality
by putting $C=B^{1/2}A^{1/2}$ and $D=B$.
\qed

We now present a new and much easier proof of the main result in \cite{ka}; moreover, we extend
its validity to include the range $0<q\le 1$ and $2\le q\le 3$.
\begin{theorem}
Let $A$ be a positive semi-definite block matrix, partitioned as above. Then for $1\le q\le 2$,
with $\beta=||B||_q$, $\delta=||D||_q$ and $\gamma=||C||_q$,
\be\label{ineq}
\trace A^q \le (2^q-2) \gamma^q + \beta^q+\delta^q.
\ee
For $0<q\le 1$ and for $2\le q\le 3$, the reversed inequality holds.
\end{theorem}
\textit{Proof.}
As already noted in \cite{ka}, it is enough to consider positive $C$. In that case
the inequality (\ref{ineq}) can be rephrased as follows:
\be\label{ineq2}
\trace\twomat{B}{C}{C}{D}^q - \trace\twomat{B}{0}{0}{D}^q
\le \trace\twomat{C}{C}{C}{C}^q - \trace\twomat{C}{0}{0}{C}^q.
\ee

Consider first the cases $0<q\le1$ and $2\le q\le 3$.
Let us calculate the minimum value of the left-hand side of (\ref{ineq2})
over all allowed $B$.
The constraint on $B$, originating from the requirement $A\ge0$, is $B\ge CD^{-1}C$.
We will show that the minimum over $B$ is obtained in $B=B_0:=CD^{-1}C$.
Let us thereto put $B=B_0+t\Delta$, with $\Delta\ge0$, and define
$$
f(t) := \trace\twomat{B_0+t\Delta}{C}{C}{D}^q - \trace\twomat{B_0+t\Delta}{0}{0}{D}^q.
$$
The Fr\'echet derivative \cite{bhatia} of $f$ is given by
$$
f'(t) = q\trace\left[\left(\twomat{B}{C}{C}{D}^{q-1}-
\twomat{B}{0}{0}{0}^{q-1}\right)\,\twomat{\Delta}{0}{0}{0}\right].
$$
Introducing the projector $P=\id\oplus 0$, we can write
$$
f'(t) = q\trace\left[\left(P\twomat{B}{C}{C}{D}^{q-1}P-\left(P\twomat{B}{C}{C}{D}P\right)^{q-1}\right)
\,\twomat{\Delta}{0}{0}{0}\right].
$$
For $0< q\le 1$ and for $2\le q\le 3$, the function $x\mapsto g(x)=x^{q-1}$ is operator convex on $(0,+\infty)$.
Therefore (\cite{bhatia}, Exercise V.2.2 applied to a compression to the upper left block)
$$
P\twomat{B}{C}{C}{D}^{q-1}P \ge \left(P\twomat{B}{C}{C}{D}P\right)^{q-1}.
$$
This shows that $f'(t)\ge0$ and that $f(t)$ is indeed minimal in $t=0$.
Therefore, we can restrict to $B=CD^{-1}C$.
The theorem now follows immediately from Corollary \ref{cor:abq3}.

For the case $1\le q\le 2$, we proceed in exactly the same way.
The function $x\mapsto g(x)=x^{q-1}$ is now operator concave on $[0,+\infty)$, with $g(0)=0$.
Therefore (\cite{bhatia}, Theorem V.2.3)
we now have
$$
P\twomat{B}{C}{C}{D}^{q-1}P \le \left(P\twomat{B}{C}{C}{D}P\right)^{q-1}.
$$
Hence, $f(t)$ is maximal in $t=0$.
\qed

We believe that the reversed inequality holds for $q>3$, but we haven't been able to prove this yet. The proof
given in \cite{ka}, using a duality argument, is incorrect.
%-----------------------------------------------------------------------
%\begin{ack}
%\end{ack}
%------------------------------------------------------------- BIBLIOGRAPHY

%%%%%%%%%%%%%%%%%%%%%%%%%%%%%%%%%%%%%%%%%%%%%%%%%%%%%%%%%%%%%%%%%%%
\end{document}